\documentclass[11pt]{article}

\title{Comments on combinatorial interpretation of Fibonomial
coefficients- an e-mail style letter (*)}
\author{A.K.Kwa\'sniewski\\  %%% ENDNAME
\\ High School of Mathematics and Applied Informatics\\
PL - 15-021 Bialystok , ul.Kamienna 17,  Poland
\\e-mail: kwandr@uwb.edu.pl}

\usepackage{amsmath,amsthm}

\chardef\bslash=`\\ % p. 424, TeXbook
\begin{document}
\maketitle \textbf{(*) Bulletin of the Institute of Combinatorics
and its Applications} vol. \textbf{42} (\textbf{2004}): 10-11. \\
Presented also at the Gian-Carlo Rota Polish Seminar\\
\noindent \emph{http://ii.uwb.edu.pl/akk/sem/sem\_rota.htm}
\vspace{0.5cm}

\textsc{I}. Up to our  knowledge -since about $126$ years we were
lacking of classical type combinatorial interpretation of
Fibonomial coefficients as it was Lukas \cite{1} - to our
knowledge -who was the first who had defined Finonomial
coefficients  and derived a recurrence for them (see Historical
Note in  \cite{2} ). Namely as accurately noticed by Knuth and
Wilf in \cite{3} the recurrent relations for Fibonomial
coefficients appeared already in $1878$ Lukas work \cite{1} and in
our opinion - Lucas  Théorie des Fonctions Numériques Simplement
Périodiques is the far more non-accidental context  for binomial
and  binomial-type coefficient - Fibonomial coefficients included.

\textsc{II}. Recently \cite{4} Kwa\'sniewski  combinatorial
interpretation of Fibonomial coefficients has been proposed in the
spirit \cite{2} of the historically classical standard
interpretations according to the schematic correspondences:

SETS ---SUBSETS (without and with repetitions)--- Binomial
coefficient ---i.e. we are dealing with LATTICE of subsets

SET PARTITIONS: Stirling numbers of the second kind
---  number of partitions into exactly $k$  blocs - i.e. we are
dealing with  LATTICE of partitions.

PERMUTATION PARTITIONS : Stirling numbers of the first kind  ---
number of permutations containing exactly  $k$  CYCLES.

SPACES: $q$-Gaussian coefficient --- number of $k$-dimensional
subspaces in n-th dimensional space over Galois field  $GF(q)$ ---
i.e. we are dealing with LATTICE of subspaces. (For nontrivial and
fruitful Konvalina`s unified interpretation of the Binomial
Coefficients, the Stirling Numbers, and the Gaussian Coefficients
see  \cite{5} ).

POSET --- Fibonomial coefficients --- number of corresponding
(see: \cite{4,2} ) finite "cobweb" subposets of the so called
"cobweb" poset.

\textsc{III}. At the time of publishing \cite{4} Kwa\'sniewski was
not aware of the existence of the relevant preprint \cite{6} of
Ira M. Gessel and   X. G. Viennot  (Just few hours ago I have
noticed this article via Google) There right after the Theorem
$25$ (see Section $10$ , page $24$  in \cite{6} ) relating  the
number N(R) of nonintersecting $k$-paths  to Fibonomial
coefficients (via $q$-weighted type counting formula) the authors
express their wish - worthy to be quoted: "\textit{it would be
nice to have a more natural interpretation then the one we have
given}"..... " \textit{R. Stanley has asked if there is a binomial
poset associated with the Fibonomial coefficients..."}

- Well. The cobweb locally finite infinite poset by Kwa\'sniewski
from \cite{2,4} is not of binomial type. Even more ; the incidence
algebra origin arguments seem to make us not to expect binomial
type poset  come into the game \cite{7}.  Am I right?

    An immediate question arises - what is the relation like between
these two: Gessel and Viennot \cite{6}  and  \cite{4} points of
view? We  shall try to elaborate more on that soon.

\end{document}